\documentclass[numbook]{svjour3}

\usepackage{amscd,amsfonts,amssymb,amsmath,latexsym,array,hhline,xcolor,graphicx}
\usepackage[utf8]{inputenc} 

\newcommand\F{\mbox{I\kern-2pt F}}

\newcommand\cP{{\cal P}}

\newcommand\R{\bbr}

\def\bbr{{\mathbb R}}

\newcommand\e{{\varepsilon}}
\newtheorem{theo}{Theorem}[section]

\newtheorem{lemm}[theo]{Lemma}
\newtheorem{coro}[theo]{Corollary}

\newtheorem{rem}[theo]{Remark}
\newcommand\fdem{$\Box$}
\newcommand\beq{\begin{equation}}
\newcommand\eeq{\end{equation}}
\newcommand\bea{\begin{eqnarray}}
\newcommand\eea{\end{eqnarray}}
\newcommand\bean{\begin{eqnarray*}}
	\newcommand\eean{\end{eqnarray*}}

\let\oldrem\rem
\let\oldendrem\endrem
\renewenvironment{rem}{\oldrem\em}{\oldendrem}
\newcommand \bma{\begin{pmatrix}}
\newcommand \ema{\end{pmatrix}}

\begin{document}
	\title{An Axiomatic Viewpoint on the Rogers--Veraart and Suzuki--Elsinger Models of Systemic Risk \\ 
	}
	\author{ Yuri Kabanov \and Arthur Sidorenko} 
	
	\institute{\at	Lomonosov Moscow State University, Institute of Informatics Problems,
	Federal Research Center “Computer Science and Control” of the Russian Academy of Sciences
	Moscow, Russia, and Laboratoire de Math\'ematiques, Universit\'e de Franche-Comt\'e, 
	16 Route de Gray, 25030 Besan\c{c}on, cedex, France \\
  \email{ykabanov@univ-fcomte.fr}. \\
	\and Lomonosov Moscow State University, Moscow, Russia\\
	 \email{Artur.Sidorenko@student.msu.ru}}

	\date{  }
	
	\maketitle
	
\begin{abstract}
	We study a model of clearing in an interbank network with cross-holdings 
and default charges.  Following the  Eisenberg--Noe approach, we  define the model via  a set of natural financial 
regulations including those related with eventual default charges and derive   a finite family of fixpoint problems. 
These  problems are parameterized by vectors of binary variables.  Our model combines features of the Ararat--Meimanjanov, Rogers--Veraart,  and Suzuki--Elsinger networks. We develop  me\-thods of computing  the maximal and  minimal clearing pairs using the mixed integer-linear programming
and  a Gaussian elimination algorithm. 

	\end{abstract}
	 		
		 

\section{Introduction}

In a financial system with strongly interconnected institutions, a shock touching even a small number of entities can propagate through the network and lead to substantial losses. Clearing is a procedure that diminishes the total amount of interbank liabilities and, as such, decreases the risk of a system-wide breakdown. In the seminal paper \cite{eisenberg2001systemic} by Eisenberg--Noe, the  notion of clearing vectors was introduced in an axiomatic way,   via  limited liability and absolute priority rules.  It was shown that the set of clearing vectors is the set of fixed points of a simple nonlinear mapping and this set contains the maximal and minimal  elements. The further development led to more sophisticated models including cross-holdings, seniority of debts, and default charges, see, e.g.  
\cite{suzuki2002valuing},  \cite{elsinger2009financial}, \cite{RV2013}, and a survey paper \cite{kabanov2018clearing}. Dynamic versions of the Eisenberg--Noe network are suggested, for instance, in \cite{feinstein2019dynamic} and \cite{djete2021mean}. 

The adding of default fees, introduced in the paper by Rogers and Veraart \cite{RV2013}, is of particular interest because it allows to study the problem of a rescue consortium to aid  insolvent institutions. 
In general, default charges have a crisis amplifying effect: external pay-offs may increase the number 
of  failures. 
In the Rogers--Veraart model the clearing vectors were defined directly as solutions of a fixpoint equation.  As was observed by Ararat and Meimanjanov in \cite{Ararat2019}, such a formulation does not coincide with a formulation defined in terms  of financial regulations. 

In this note we consider  a model with default charges and cross-holdings using the Ararat--Meimanjanov
``axiomatic" approach. Our aim is to derive the fixpoint equations and suggest algorithms to find their minimal and maximal solutions.  In particular, we consider a method using the mixed integer-linear programming
and  a Gaussian-type dimensionality reduction algorithm. 

\smallskip {\bf Notations:} For vectors $a, b \in \bbr^N$ the symbol $a\le b$ denotes the component-wise partial ordering $a^i \leq b^i$ for every $i = 1, \dots, N$,  $[a,b]:=\{x\in \bbr^N\colon a\le x\le b\}$, ${\bf 1}_{\{a<b\}}$
is the vector with the components  $I_{\{a^i< b^i\}}$,  ${\bf 1}:=(1,\dots,1)$, ${\bf 0}:=(0,\dots,0)$.    
The symbol  $a\circ b:=(a^1b^1,\dots,a^Nb^N)$ stands for the Hadamard (component-wise) product,

\section{Suzuki--Elsinger model with default payments}

Let us consider a financial network consisting of $N>1$ banks. The bank $i$ has a cash reserve  $e^i$, a liability $l^{ij}$ towards the bank $j\neq i$ with the total  $l^i=\sum_j l^{ij}$, and possess a share $\theta ^{ji}$ of the bank $j$.  
It is convenient to introduce the   relative liabilities matrix $\Pi=(\pi^{ij})$ with 
\beq
\label{sm}
\pi^{ij}:=\frac {l^{ij}}{l^i}=\frac {l^{ij}}{\sum_j l^{ij}}, \quad\hbox{if \ }l^i\neq 0, \ \hbox{and}\  \pi^{ij}:=\delta^{ij}\  \hbox{otherwise}, 
\eeq
where the Kronecker symbol $\delta^{ij}=0$ for $i\neq j$ and $\delta^{ii}=1$.
Then $\pi^{ij}$  describes the fraction of the value of the debtor $i$ due to the creditor $j$ of the total interbank debt of $i$.  The value $\pi^{ii}=1$ means that 
the bank $i$ has no interbank debts. 

The matrix $\Theta=(\theta^{ij})$ is assumed to be substochastic and such that unit is not its eigenvalue. 

 Fix $\alpha, \beta, \gamma \in [0,1]$.  
It is assumed that there is  Central Clearing Counterparty (CCP) calculating  {\it clearing} and {\it equity}  vectors, 
denoting $p$ and $V$, and provide  the settlement service for creditors and debtors.   

A {\it clearing pai}r $(p,V)\in [0,l]\times  \bbr^N_+$ is determined by the following rules:    

\begin{enumerate}
    \item    
    {\bf Limited liability}: $p^i \leq (e + \Pi' p + \Theta' V)^i$ for every $i$,
    \item
    {\bf Absolute priority}: either $p^i = l^i$ or $p^i = (\alpha e + \beta \Pi' p + \gamma \Theta' V)^i$. 
    \item
    {\bf Equity evaluation}: if $p^i = l^i$, then $V^i = (e + \Pi' p + \Theta' V - p)^i$, otherwise $V^i = 0$.
\end{enumerate}

The  first rule means that the clearing payment cannot exceed the available resources 
(cash $e^i$ plus collected debts $\sum_j \pi^{ji}p^j$ plus the total of owned shares $\sum_j \theta^{ji} V^j$) while
the second rule stipulates that either debts are payed in full, or all resources are distributed with the charges payed in the case of the default, namely, the amount 
$$
(1-\alpha)e^i + (1-\beta) (\Pi'p)^i+(1- \gamma)( \Theta' V)^i.
$$ 
The third relation is just the definition of the equity:  $V^i$ is the total value of  assets of the bank $i$ after clearing  (corresponding to the clearing vector $p$). 
By virtue of the limited liability rule, all components $V^i\ge 0$. 

We denote by $\cP$ the set of clearing pairs.

 Put 
\bean
G(x, V) &:= & (e + \Pi' x + \Theta' V - l) \circ {\bf 1}_{\{x = l\}},\\
G_+(x, V) &:= & (e + \Pi' x + \Theta' V - l)^+ \circ {\bf 1}_{\{x = l\}}. 
\eean
In this notation, if $(p,V)$ is  a clearing pair, (i.e. a clearing vector and an equity vector), then   $V = G(p, V)$.

\begin{lemm}
\label{lem21}
If  $x \in \R^N$, then  the equations $V = G(x, V)$ and $V = G_+(x, V)$ have unique solutions.  As functions of $x$, the solution of the first equation is  affine, the solution of the second is  convex increasing. 
\end{lemm}
\noindent 
{\sl Proof.}
Put $\Lambda := {\rm diag}\;{\bf 1}_{\{x=l\}}$, $a(x): = \Lambda (e + \Pi'x - x)$, and $B := \Theta \Lambda$.
In this notation
$$ 
G(x, V) = a(x) + B' V. 
$$
The  linear equation $V=G(x,V)$ has a unique solution if the matrix 
$I-B$ is invertible, i.e. unit is not an eigenvalue of $B$ (recall that this property is  assumed).

 Analogously, $G_+(x, V)$ can be written in the form 
$$ G_+(x, V) = (c(x) + B' V)^+, $$
for $c(x) := \Lambda (e + \Pi'x - l)$. 
It remains to apply   Lemma 4.3 from \cite{kabanov2018clearing} claiming that for every $y\in \R^N$ the equation $v=(y+\tilde \Theta v)^+$, where $\tilde \Theta$ is a substochastic matrix with invertible $I-\tilde \Theta$, has a unique solution $v=v(y)$ which is a monotone increasing  convex function and so is the function $H(x):=v(c(x))$.    \fdem 

\smallskip
\noindent 
{\sl Remark.}
In the above lemma, the existence can be deduced from the Knaster--Tarski theorem. Indeed, put 
$z := e + \Pi' l - l$, $K := (I - \Theta')^{-1} z^+$. Then for any $p \in [0, l]$
$$
G_+(p, K) \leq G_+(l, K) = (z + \Theta' K)^+ \leq z^+ + \Theta' K =  K.  
$$
It follows that  the monotone function  $V \mapsto G_+(p, V)$ maps $[0, K]$ into itself.

\begin{lemm}
\label{lem22}
If $(p,V)\in\cP$, then  
$V = G_+(p, V)$.
\end{lemm}
\noindent
{\sl Proof.}
As we already observed, the equity evaluation means that   $V = G(p, V)$  and     $p \le e + \Pi' p + \Theta' V$ because of the limited liability. The components  $G^i(p, V)$ and  $G^i_+(p, V)$
coincide if $p^i=l^i$, and equal to zero of $p^i<l^i$. 
\fdem

\smallskip

If $\alpha=\beta=1$ and the  matrix $\Theta$ is zero (no crosshodings), the model is reduced to that of  Eisenberg--Noe.  If $\alpha,\beta\in ]0,1]$ and  $\Theta$ is zero,  we get the Ararat--Meimanjanov 
model. The model of Suzuki--Elsinger has a non-trivial substochastic matrix $\Theta$ and  $\alpha=\beta=\gamma=1$, i.e.  default charges are not imposed. The relation with the Rogers--Veraart model we discuss later. 

\smallskip

To determine clearing pairs we  use the absolute priority rule. For $p\in [0,l]$
we put   
\bean
\Phi^i_0(p, V) &:= & (\alpha e+\beta \Pi' p + \gamma \Theta' V)^i \wedge l^i,\\
\Phi^i_1(p, V)& := & (\alpha e+\beta \Pi' p + \gamma \Theta' V)^i d^i+ (1 - d^i)l^i,
\eean
where 
$$
d^i := I_{\{ (e + \Pi' p + \Theta' V)^i < l^i \}}, \qquad i=1,\dots,N. 
$$ 

For a binary vector  $b=(b^1,...,b^N) \in \{0, 1\}^N$ we define on $\R^N\times \R^N$ the $\R^N$-valued function  $(p,V)\mapsto \Phi_b (p, V)$ with 
$$
\Phi_b (p, V) = (\Phi^1_{b^1}(p, V), \dots, \Phi^N_{b^N}(p, V) ).
$$ 

In particular, for ${\bf 0}:=(0,...,0)$ and ${\bf 1}:=(1,...,1)$ we have 
\bean 
\Phi_{\bf 0}(p, V) &:=&  (\alpha e+\beta \Pi' p  + \gamma \Theta' V)\wedge l, \\ 
\Phi_{\bf 1}(p, V) &:=& (\alpha e+\beta \Pi' p + \gamma \Theta' V)\circ d  + ({\bf 1} - d) \circ l,
\eean 
where $d = (d^1, \dots, d^N)$. 

\smallskip
Let  $H(x)$ be the fixpoint $G_+(x, \cdot)$. By  Lemma \ref{lem22} the mapping $H: \bbr^N \to \bbr^N_+$ is monotone.
Let $F^i_j(p) := \Phi^i_j(p, H(p))$, $i=1,\dots, N$,   $j=0,1$, and let  
$$
F_b(p) := \Phi_b(p, H(p))=( \Phi^1_b(p, H(p)),\dots,  \Phi^N_b(p, H(p))), \qquad b\in \{0,1\}^N.  
$$
It is easily seen that  $(b,p)\mapsto F_b(p)$ is a monotone mapping of $\{0,1\}^N\times [0,l]$ into $[0,l]$ with respect to the componentwise ordering.
By the Knaster--Tarski theorem, see e.g. \cite{kabanov2018clearing}, the set of its fixpoints is non-empty and contains the minimal and maximal elements. By the same reason, for every $b\in \{0,1\}
$ the set of  fixpoints of the mapping $p\mapsto F_b(p)$ is non-empty and contains the minimal and maximal elements.
\smallskip

Let  $S:=\{p \in [0, l]\colon (p,V)\in \cP\}$, i.e. $S$ is the projection of the subset $\cP$ of the ``plane'' $\R^N\times \R^N$ into the ``$x$-axis'' 
$\R^N$. 

The following  result gives a characterization of the set $S$ of  clearing vectors. 

\begin{theo}
\label{RV_t1}
 $S=\bigcup_b\{p \in [0, l]\colon p=F_b(p)\}$.
\end{theo}
\noindent
{\sl Proof.}
\smallskip
Define the set $
 W_b := \bigcap_{i=1}^{N} W^i_{b^i}$, $b \in \{0, 1\}^N$, where  $W^i_j := \{ p \in S \colon p^i = F^i_j (p)\}$. 
 
 Let $p\in S$, that is $(p,V) \in \cP$ for some $V=H(p)$,  see Lemma \ref{lem22}. There are two cases  $p^i = l^i$ and $p^i < l^i$. In the first one 
 $$l^i \le (e + \Pi' p + \Theta' V)^i=(e + \Pi' p + \Theta' H(p))^i,
 $$
implying that that $F^i_1 (p)=l^i=p^i$, hence, $p \in W^i_1$.
In the second case,   by the absolute  priority rule, $p^i = (\alpha e + \beta \Pi' p + \gamma \Theta' H(p))^i$, that is, $p^i = F^i_0(p)$ and $p \in W^i_0$. Thus, 
 $S = W^i_0 \cup W^i_1$ whatever is $i$ and  
 $$
 S=\bigcap_i (W^i_0 \cup W^i_1)=\bigcup_{b \in \{0, 1\}^N}  W_b\subseteq   \bigcup_b\{p \in [0, l]\colon p=F_b(p)\} .
$$

\smallskip

To prove the converse inclusion,  fix some $b$ and take $x \in [0, l]$ such that $x = F_b(x)$. Put $V = H(x)$. It is easily seen that $(x, V)$ satisfies all the axioms.
\fdem

\begin{coro}
\label{RV_c1}
Let  $\underline p_b$ and $\bar p_b$ be the minimal and maximal  fixpoints of the mapping $F_b$. Then 
 $\underline p_{\bf 0}$ and $ \bar p_{\bf 1}$ are the minimal and maximal  fixpoints of $S$.
 \end{coro}
\smallskip
\noindent
{\sl Proof.} 
Since functions $F_b$ are increasing in $b$, so are the $\underline p_b$ and $\bar p_b$, see, e.g., \cite{kabanov2018clearing}. This implies the result. \fdem

\smallskip

Let us elaborate on the above statements. The maximal clearing vector $p^*$ is the  component of the vector  $(p^*, V^*)\in [0,l]\times \bbr^N_+$ which is the maximal solution of the system 
\bean
    p& =&  (\alpha e+\beta \Pi' p + \gamma \Theta' V) \circ d + l \circ ({\bf 1} - d), \\
     V &= & ( e + \Pi' p +  \Theta' V - l)^+ \circ ({\bf 1} - d),
\eean
where $d := {\bf 1}_{\{e + \Pi' p + \Theta' V < l \}}$. 

The minimal clearing vector $p_*$ is the  component of the vector  $(p_*, V_*)\in [0,l]\times \bbr^N_+$ which is the minimal solution of the system 
\bean
    p & =&  (\alpha e+\beta \Pi' p + \gamma \Theta' V) \wedge l, \\
     V & =&  ( e + \Pi' p +  \Theta' V - l)^+ \circ {\bf 1}_{\{p = l\}} .
\eean

Consider two particular cases. Let $\Theta = 0$. Then  $V$ is irrelevant. The equation for the maximal fixpoint is 
$$
p =   (\alpha e+\beta \Pi' p ) \circ d + l \circ ({\bf 1} - d),
$$
where $d := {\bf 1}_{\{e + \Pi' p < l \}}$. 

This equation coincides with that introduced in \cite{RV2013} as the definition of the clearing vector. 
It means that the maximal clearing vector in the ``axiomatic'' definition of \cite{Ararat2019} and in the  definition via fixpoint of \cite{RV2013} coincide. However, in the axiomatic approach, the minimal clearing vector is the minimal solution of a different equation
$$ 
p = (\alpha e+\beta \Pi' p ) \wedge l
$$
and may be different  from that coming via fixpoint equation of \cite{RV2013}. 

Our arguments are intended to show that  the  ``axiomatic'' description also leads 
to a fixpoint problem allowing, in principle, to find all clearing vectors and corresponding equities. Unfortunately, the number of equations are exponentially 
growing. E.g. for a small financial system with only ten banks  we have  20480 equations. On the other hand, it seems that the largest clearing vector is of a major practical interest and to get it one  can consider the system of only $2N$ equations suggested by Rogers--Veraart.

Comparably to the  Rogers--Veraart approach, the axiomatic one might seem rather perplexing because it leads to $2^N$ equations instead of just one. There is no  hope that a significant proportion of these equations is redundant and one can rule out them. Indeed, take arbitrary $\Pi$, $l$ with all  $l^i > 0$, put $\Theta = 0$, $e = l$, $\gamma = 1$, $\beta=\alpha$. Choose $\alpha$ such that that $\alpha (l + \Pi' l)^i < l^i$ for all $i$. Obviously, the solution of the equation $p = F_b(p)$ has  the component  $p^i = l^i$ if $b^i=1$ and $p^i < l^i$, otherwise. In this example the sets $W_b$ are disjoint singletons.

Now let us consider the case where $\alpha = \beta = \gamma = 1$ as in the  Suzuki--Elsinger model. In this case, $F_{\bf 0} = F_{\bf 1}$ and clearing vectors are solutions of the system 
\bean
 p &=& (e + \Pi' p + \Theta' V) \wedge l, \\
V &=& (e + \Pi' p + \Theta' V - l)^+ \circ {\bf 1}_{\{p = l\}}.
\eean
Obviously, the second equation can be replaced by $V = (e + \Pi' p + \Theta' V - l)^+$. Hence, the proposed ``axiomatic" definition coincides with that of Suzuki--Elsinger model in \cite{elsinger2009financial}.

\smallskip
A short remark on the interpretation of the parameter $b$.  
Assume for simplicity that $\alpha = \beta = \gamma < 1$.    Let   $(p, V)$ be a clearing pair corresponding to some binary vector $b$ and let  $x(p,V) := e + \Pi' p + \Theta'V$. The component  $x^i(p,V)$ is the total  assets of bank $i$. For $z \in \bbr_+$, put $h^1_0(z) := (\alpha z) \wedge l^1$ and $h^1_1 (z) := \alpha z {\bf 1}_{ \{ z < l^1 \} } + l^1 {\bf 1}_{ \{ z \geq l^1 \} }$.The value $h^1_{b^1}(x^1(p, V)) = \Phi^1_{b^1}(p, V)$  is the debt payment of bank 1.  If $x^1(p, V) \in [l^1, l^1/\alpha)$, then  $h^1_0(x^1(p, V)) < h^1_1(x^1(p, V)) = l^1$. 
Let us compare clearing pairs $(p,V)$ and $(\bar p,\bar V)$ corresponding to the binary vectors $(0,\tilde b)$ and $(1,\tilde b)$ with all components equal except the first one. Then $(p,V)\le (\bar p,\bar V)$ and 
$x(p,V)\le x(\bar p,\bar V)$.  It may happen that 
$ x^1(\bar p, \bar V) \in [l^1, l^1/\alpha)$ but $x^1( p, V) <l^1$.  This leads to the conclusion that the banks are motivated to require  CCP to calculate maximal  clearing pairs using $b={\bf 1}$.

\section{Clearing pairs via integer-linear programming }

In this section we inroduce two integer-linear programming problems to find the maximal and  and minimal clearing  clearing pairs in the models with crossholdings and default charges. To this end we consider the linear function $f: \{ 0, 1 \}^N \times [0, l] \times \bbr^N_+ \to \bbr$ with 
$$
f(a, p, V) = \sum_{i=1}^N f_{1i} a^i + f_{2i} p^i + f_{3i} V^i,
$$
where all coefficients  $ f_{1i}$,  $f_{2i}$, and $f_{3i}$ are strictly positive. Put 
$$
x = x(p, V) = e + \Pi' p + \Theta' V, \qquad y = y(p, V) = \alpha e + \beta \Pi' p + \gamma \Theta' V.
$$

\subsection{The maximal clearing pair}
Let $\kappa:=||K||_\infty$ where $K$ is the vector introduced in the remark after Lemma \ref{lem21}. 
 Then $H(l) \leq \kappa {\bf 1} $ where the function $H$ was defined in Section 2. 

Problem {\bf P1}: maximize $f$ under the constraints 
\bea
 \label{max1}
 &&p \le y(p,V) + a \circ l, \\
 \label{max2}
&& a \circ l \le x(p,V), \\
\label{max3}
&& V \le x(p,V) - a \circ  l, \\
\label{max4}
&& V \leq \kappa a.
\eea
Immediately, if $(p, V)$ is a clearing pair, then $({\bf 1}_{\{p = l\}}, p, V)$ is admissible for {\bf P1}.

\begin{lemm}
If $(\hat a, \hat p, \hat V)$ solves {\bf P1}, then $\hat a = {\bf 1}_{\{\hat p = l\}}$.
\end{lemm}
{\sl Proof.}  Suppose that 
 $\hat p^j < l^j$ but $\hat a^j=1$. In this case, we have from (3.1) that   $p^j\le y^j(\hat p,\hat V) + l^j $. We can replace  $\hat p^j$ by a larger value $\hat p^j+\e$ with $\e\in (0,l^j-\hat p^j)$ without violating the constraints.  Since the cost function $f$ is strictly increasing we get a contradiction with
 the optimality of  $(\hat a, \hat p, \hat V)$. Thus,  $\hat a^j = 0$. 

Suppose that  $p^j = l^j$ but $\hat a^j = 0$. Then the component $\hat V^j = 0$ by virtue of  \eqref{max4} and  $l^j \leq y^j (\hat p,\hat V)\leq x^j(\hat p, \hat V)$ due to   \eqref{max1}. Put $\tilde a^j = 1$, $\tilde a^i = \hat a^i$ for $i \neq j$. Then  $(\tilde a, \hat p, \hat V)$ is admissible and we again obtain a contradiction with the optimality. \qed

\begin{theo}
If $(\hat a, \hat p, \hat V)$ solves {\bf P1}, then $(\hat p, \hat V)$ is the maximal clearing pair.
\end{theo}
{\sl Proof.} 
Let $(\hat a, \hat p, \hat V)$ solve {\bf P1}. 
It suffices to verify that $(\hat p, \hat V)$ is a clearing pair.
 Due to the above lemma we already know that $\hat a = {\bf 1}_{\{\hat p = l\}}$. Let us check that the required rules are met. 

\smallskip
 {\bf Limited liability}. If $\hat p^j < l^j$, then $\hat p^j \leq y^j(\hat  p,\hat  V)$ from \eqref{max1}. Thus, $\hat p^j \leq x^j(\hat  p,\hat  V)$. If $\hat p^j = l^j$, then $\hat a^j = 1$ and  $l^j \leq  x^j(\hat  p,\hat  V)$ from \eqref{max2}.

\smallskip
 {\bf Absolute priority}. If $\hat p^j = l^j$, then there is nothing to prove. If $\hat p^j < l^j$, then $\hat p^j \leq y^j(\hat  p,\hat  V)$ by \eqref{max1}. If the inequality were strict, one could add a small $\e>0$ to $\hat p^j$ and obtain again an admissible vector.

\smallskip
 {\bf Equity evaluation}. If $\hat p^j < l^j$ and $\hat a^j = 0$, then $\hat V^j = 0$ because of \eqref{max4}. If $\hat p^j = l^j$ and $\hat a^j = 1$, then  $l^j \leq x^j(\hat  p,\hat  V)$ and $\hat V^j \leq x^j(\hat  p,\hat  V) - l^j$ from \eqref{max2}. In the case where the latter inequality is strict, one can add a small $\e>0$ to $\hat V^j$ and get an admissible vector. \fdem

\subsection{The minimal clearing pair}

Put  $\kappa_1 := || e + \Pi'l +\kappa \Theta' {\bf 1} ||_\infty $. Then $e + \Pi'l + \Theta' H(l) \leq \kappa_1 {\bf 1}$.  
Problem {\bf P2}: minimize $f$ under the constraints 
\bea
 \label{min1}
&& p \geq y(p,V) - \kappa_1 a, \\
 \label{min2}
&& p \geq a \circ l, \\
\label{min3}
&&  ( {\bf 1} - a) \circ l  + \kappa_1 a \geq y(p,V) , \\
\label{min4}
&& V \geq x(p,V) - l - \kappa({\bf 1} - a).
\eea

 We will show, under mild assumptions, that the solution $(\check a, \check p, \check V)$ of {\bf P2} is the minimal clearing pair. More precisely, it is the minimal solution of the system 
\bea
    \label{min_eq1}
    p &=& \Phi_{\bf 0} (p, V), \\
    \label{min_eq2}
    V &=& \check G_+(p, V),
\eea
where 
$$
\check G_+ (p, V) = (x(p,V) - l)^+ \circ {\bf 1}_{\{y(p,V)> l\}}.
$$

Note that  for every solution $( p,  V)$ of the above system, $({\bf 1}_{\{y(p,V)> l\}}, p, V)$ is admissible for {\bf P2}.

\begin{lemm}
The vector $\check p = ({\bf 1} - \check a) \circ x(\check p, \check V) +  \check a \circ l$.
\end{lemm}
\noindent
{\sl Proof.} If $\check a^j = 0$, then by \eqref{min1} $\check p^j \geq y^j(\check p,\check V)$. If the inequality is strict, one can decrease $\check p^j$ by a small $\e>0$. Hence, 
$\check p^j = y^j(\check p,\check  V)$. If
 $\check a^j = 1$, then by \eqref{min2}, $\check p^j = l^j$.  \qed

\begin{lemm}
The vector  $\check V = \check a \circ (x(\check  p,\check  V)-l)^+$.
\end{lemm}
{\sl Proof.}  Let $\check a^j = 0$. Then $\check V^j = 0$ (otherwise  $(\check a, \check p, \check V)$ cannot be the solution of  {\bf P2}). 

Let $\check a^j = 1$ and $x^j(\check  p,\check  V) - l^j \geq 0$. From \eqref{min4}, $\check V^j \geq x^j(\check  p,\check  V) - l^j$. Since $f$ attains its minimum at $(\check a, \check p,\check  V)$, necessarily, $V^j = x^j(\check  p,\check  V) - l^j$. 

Let $a^j = 1$, but $x^j(\check  p,\check  V) - l^j < 0$. Then $\check V^j = 0$ by the same reasoning.   \qed

\begin{lemm}
The vector $\check a = {\bf 1}_{\{y(\check  p,\check  V) > l \} }$. In particular, $\check p = y(\check  p,\check  V) \wedge l$.
\end{lemm}
{\sl Proof.} 
Let $y^j (\check  p,\check  V)\leq l^j$.  Assume for a minute that $\check a^j = 1$. Immediately, $\check p^j = l^j$. Put $\tilde a^j = 0$ and $\tilde a^i = \check a^i$ for $i \neq j$, then $(\tilde a, \check p, \check V)$ is admissible. A contradiction. \\
Let  $y^j(\check  p,\check  V) > l^j$. Then $\check a^j = 1$ by \eqref{min1}. \qed

\begin{theo}
Let $(\check a, \check p, \check V)$ be the solution of {\bf P2}. Then $(\check p, \check V)$ is the solution of the system \eqref{min_eq1} -- \eqref{min_eq2}. In particular, $\check p \leq \underline p_{\bf 0}$, where $\underline p_{\bf 0}$ was defined in Corollary \ref{RV_c1}. If $y^j(\check p, \check V) \neq l^j$ for every $j$, then $\check p =\underline p_{\bf 0}$.
\end{theo}
{\sl Proof.}
The fact that $(\check p, \check V)$ solves \eqref{min_eq1} -- \eqref{min_eq2} follows from the above lemmata. Note that $(\underline p_{\bf 0}, H(\underline p_{\bf 0}))$ solves the system
\bean
p &=& y(p, V) \wedge l, \\
V &=& (x(p, V) - l)^+ \circ {\bf 1}_{\{x=l\}}
\eean
which  is equivalent to the following one:
\bea
\label{min_eq3}
p &=& y(p, V) \wedge l, \\
\label{min_eq4}
V &=& (x(p, V) - l)^+ \circ {\bf 1}_{\{y(p, V) \geq l\}}.
\eea
Put $\tilde G_+(p, V) := (x(p, V) - l)^+ \circ {\bf 1}_{\{y(p, V) \geq l\}}$. Obviously, $\check G_+(p, V) \leq \tilde  G_+(p, V)$ for every $(p, V)$. Due to the Knaster--Tarski theorem, we have $\check p \leq \underline p_0$.
Furthermore, if $y^j(\check p, \check V) \neq l^j$ for every $j$, then $\check G_+(p, V) = \tilde  G_+(p, V)$. Under this condition the pair $(\check p, \check V)$ is the solution of  \eqref{min_eq3} -- \eqref{min_eq4}, that is why $\check p = \underline p_{\bf 0}$. \fdem

\section{Gaussian method}

In this section we present a Gaussian-type   algorithm to find the maximal clearing vector in a model with crossholding and 
default charges assuming that  $\alpha=\beta=\gamma$ and $||\Theta||_\infty<1$, 
 i.e. $\sum_j \theta^{ij} < 1$ for every $i$.

In such a case   the mappings $\Phi_{\bf 1}:[0, l]\times  \R^N_+\to [0,l]$ and $G_+:[0, l]\times  \R^N_+\to \R^N_+$ are given by the formulae 
\bean
 \Phi_{\bf 1}(p, V)&: =& \alpha (e +  \Pi' p +  \Theta' V) \circ {\bf 1}_{D} + l \circ {\bf 1}_{\bar D}, \\
G_+(p, V) &:= &(e + \Pi' p + \Theta' V - l)^+ \circ {\bf 1}_{\bar D},
\eean
where $D := \{ e + \Pi' p + \Theta' V <  l  \}$. 

Our aim is to find the maximal solution of the  system 
\bea
\label{eq1}
p=\Phi_{\bf 1}(p, V) , \\
\label{eq2}
V=G_+(p, V). 
\eea
existing in virtue of the Knaster--Tarski theorem. 

In the case where 
 $(e + \Pi' l + \Theta' H(l))^i \ge   l^i $ for all $i$,    the maximal clearing pair  is 
$(l,H(l))$. Of course, to check the condition  one needs to compute the value $H(l)$ which is the unique solution of the system  $v=(c(l)+\tilde \Theta v)^+$ with $c(l):=e+P'l-l$.  It is known that the latter system
can be solved in a finite number of steps, e.g.,  by a   Gaussian-type algorithms,  \cite{kabanov2018clearing}. 

In the opposite case where at least one inequality fails to be true,  we  assume, without loss of generality, that  $e^1 + (\Pi' l)^1 + (\Theta' V(l))^1 < l^1$ and $V^1=0$. The first equation  becomes linear:  
\beq
\label{s3_eq_for_p1}
\alpha e^1 + \alpha (\Pi' p)^1 + \alpha (\Theta' V)^1 = p^1.
\eeq
 The vector equation $p=\Phi_{\bf 1}(p, V)$ can be written as the system 
\bean
p^1&=&\alpha e^1+ \alpha  \pi^{11}p^1+  \alpha  T' \tilde p +  \alpha  M' \tilde V,\\
\label{tildep}
\tilde p&=&\alpha(\tilde e+p^1R'+\tilde \Pi'\tilde p +\tilde \Theta'\tilde V)\circ \tilde 1_{D} + \tilde l \circ \tilde 1_{\bar D},
\eean
where $\tilde \Pi := (\pi^{ij})_{i, j = 2}^N$, $R: = (\pi^{1j})_{j=2}^N$, $T: = (\pi^{i1})_{i=2}^N$, $\tilde \Theta: = (\theta^{ij})_{i, j = 2}^N$, $M := (\theta^{i1})_{i=2}^N$,  and the tilde indicates  vectors with the removed first component. 

If $\alpha \pi^{11}\neq 1$, then 
\beq
\label{p1}
p^1 = \alpha(1 - \alpha \pi^{11})^{-1} (e^1 +   T' \tilde p +  M' \tilde V),
\eeq
Substituting this expression, we get that  
$$
\tilde e+p^1R'+\tilde \Pi'\tilde p +\tilde \Theta'\tilde V=e_1 +  \Pi'_1 \tilde p + \Theta'_1 \tilde V,
$$
where 
\bean
\Pi_1 &:=& \tilde \Pi + \alpha(1 - \alpha \pi^{11})^{-1}  TR, \\
\Theta_1 &:=& \tilde \Theta + \alpha(1 - \alpha \pi^{11})^{-1}  MR, \\ 
e_1 &:=& \tilde e + \alpha(1 - \alpha \pi^{11})^{-1}  e^1 R'.
\eean
It follows that $(\tilde p,\tilde V)$ is the solution of the system  
\bean
\tilde p&=&\alpha(e_1 +  \Pi'_1 \tilde p + \Theta'_1 \tilde V) \circ \tilde 1_{D} + \tilde l \circ \tilde 1_{\bar D}, \\
\tilde V&= &(e_1 + \Pi'_1 \tilde p + \Theta'_1 \tilde V - \tilde l)^+ \circ \tilde 1_{\bar D}.  
\eean

Note that $TR \tilde {\bf 1} = T (R \tilde {\bf 1}) = ((\Pi{\bf 1})^1 - \pi^{11}) T$  and, therefore,   
$$
\Pi_1\tilde {\bf 1}=\tilde \Pi \tilde {\bf 1} + (\alpha (\Pi\tilde {\bf 1})^1 - \alpha \pi^{11}) (1 -\alpha \pi^{11})^{-1} T \le \tilde \Pi \tilde {\bf 1} + T =\widetilde {\Pi{\bf 1}}\le \tilde {\bf 1}. 
$$
Thus, the matrix $\Pi_1$ is substochastic. By the same arguments we get that $||\Theta_1||_\infty<1$.
Moreover,  $ \tilde {\bf 1}_{D}={\bf 1}_{D_1}$, where $
D_1 := \{ e_1 + \Pi'_1\tilde  p + \Theta'_1 \tilde V <  l  \}$.

If $\alpha=1$ and $\pi^{11}=1$, then we can take $p^1=l^1$ and continue with the system  
\bean
\tilde p&=&(e_1 +  \Pi'_1 \tilde p + \Theta'_1 \tilde V) \circ \tilde 1_{D} + \tilde l \circ \tilde 1_{\bar D}, \\
\tilde V&= &(e_1 + \Pi'_1 \tilde p + \Theta'_1 \tilde V - \tilde l)^+ \circ \tilde 1_{\bar D},  
\eean 
where $e_1:=\tilde e+e^1R'$, $ \Pi_1:=\tilde P$, $\Theta_1:=\tilde \Theta$.

In both cases  we reduced the $2N$-dimensional problem to find the maximal solution  to $2(N-1)$-dimensional 
one of the same type. 
\bigskip
 
{\bf Acknowledgements.} The research is funded by the grant of RSF $n^\circ$ 20-68-47030 ''Econometric and probabilistic methods for the analysis of financial markets with complex structure``.


\end{document}